\documentclass[12pt]{amsart}
\usepackage{amsfonts} 
\usepackage{amssymb} 
\usepackage{verbatim} 

\def\<#1>{\langle#1\rangle} 
\def\(#1){\left(#1\right)} 
\def\forces{\Vdash} 
\def\proves{\vdash} 
\def\iff{\leftrightarrow} 
\newcommand{\from}{\mathbin{\vbox{\baselineskip=2pt\lineskiplimit=0pt
                         \hbox{.}\hbox{.}\hbox{.}}}} 

\def\corners#1{\ulcorner{#1}\urcorner} 

\newcommand{\cof}{\mathop{\rm cof}} 
\newcommand{\dom}{\mathop{\rm dom}} 
\newcommand{\ran}{\mathop{\rm ran}} 
\newcommand{\tc}{\mathop{\rm tc}} 
\newcommand{\rank}{\mathop{\rm rank}} 

\renewcommand{\P}{{\mathbb P}} 
\newcommand{\Q}{{\mathbb Q}} 
\newcommand{\M}{{\mathcal M}} 

\newcommand{\GA}{\hbox{\sc ga}} 
\newcommand{\BA}{\hbox{\sc ba}} 
\newcommand{\GCH}{\hbox{\sc gch}} 
\newcommand{\AC}{\hbox{\sc ac}} 
\newcommand{\ORD}{\hbox{\sc ord}} 
\newcommand{\ZFC}{\hbox{\sc zfc}} 
\newcommand{\ZC}{\hbox{\sc zc}} 

\newcommand{\HOD}{\hbox{\sc hod}} 
\newcommand{\VHOD}{\hbox{\sc v=hod}} 
\newcommand{\nVHOD}{\hbox{\sc v$\neq$hod}} 
\newcommand{\CCA}{\hbox{\sc cca}} 
\newcommand{\PFA}{\hbox{\sc pfa}} 

\newtheorem{Thm}{Theorem}
\newtheorem{Cor}[Thm]{Corollary}
 
\newtheorem{Sublemma}{Lemma}[Thm] 
 
\numberwithin{Claim}{Thm} 
\newtheorem{Def}[Thm]{Definition}
\newtheorem{Observation}[Thm]{Observation}
\newtheorem{TestQuestion}[Thm]{Test Questions}
\newtheorem{Question}[Thm]{Question}

\begin{document}

\title{The Ground Axiom}
\author{Jonas Reitz}
\date{July 2006}
\address{J. Reitz, The New York City College of Technology, Mathematics, 300 Jay Street, Brooklyn, NY 11201}
\email{jonasreitz@gmail.com}

\begin{abstract}
A new axiom is proposed, the Ground Axiom, asserting that the universe is not a nontrivial set forcing extension of any inner model.  The Ground Axiom is first-order expressible, and any model of $\ZFC$ has a class forcing extension which satisfies it.  The Ground Axiom is independent of many well-known set-theoretic assertions including the Generalized Continuum Hypothesis, the assertion $\VHOD$ that every set is ordinal definable, and the existence of measurable and supercompact cardinals.  The related Bedrock Axiom, asserting that the universe is a set forcing extension of a model satisfying the Ground Axiom, is also first-order expressible, and its negation is consistent.
\end{abstract}

\maketitle

\section{Introduction}

Forty years of forcing have illustrated its efficiency and flexibility in producing models of set theory with a wide variety of properties. Each such example further enriches our knowledge of the collection of models obtainable by forcing.  This collection possesses such diversity of models and such intricate structure between them that its exploration will doubtlessly continue into the forseeable future.  

I am interested in exploring the boundaries and limitations of this collection, by considering under what circumstances is the universe \emph{not} a forcing extension.  This motivates the following axiom, formulated jointly with Joel Hamkins:

\begin{Def}\rm\label{D:GAformal}The \emph{Ground Axiom} ($\GA$) is the assertion that the universe of sets $V$ is not a forcing extension of any inner model $W$ by nontrivial forcing $\P\in W$.
\end{Def}

Some observations are in order.  As stated, the Ground Axiom is second-order in nature, requiring quantification over classes $W$.  However,  I will show in Section \ref{S:GAExpressible} that the Ground Axiom has a first-order equivalent.  In addition, the Ground Axiom refers only to \emph{nontrivial} forcing, else any model is a forcing extension of itself, and only to \emph{set} forcing, a restriction that will be important in the first-order expressibility results.

There are a number of well known models of the Ground Axiom, including many of the `canonical models.'    

\begin{Observation}\label{T:CanonicalModelsofGA} The constructible universe $L$, the model $L[0^\#]$, and the canonical model of a measurable cardinal $L[\mu]$ all satisfy the Ground Axiom.
\end{Observation}

\begin{proof}In each case, the result follows from the uniqueness and minimality properties of the model.  If $L=W[g]$ were a forcing extension of an inner model $W$, then by the absoluteness of $L$ we would have $L\subset W$, and so $L=W$.  For $L[0^\#]$, suppose $L[0^\#]=W[g]$.  It is well known that $0^\#$ cannot be created by set forcing (for example, this appears in Jech \cite{Jech:SetTheory3rdEdition} as exercise 18.2, page 336), and so $0^\# \in W$.  Thus $L[0^\#]\subset W$, and so $W=L[0^\#]$.

In the case of $L[\mu]$, suppose $L[\mu]=W[g]$, where $g$ is $W$-generic for a poset $Q\in W$.  Fix a $Q$-name $\dot{\mu}$ such that $(\dot{\mu})_g = \mu$, and $q\in Q$ such that $ q \forces \dot{\mu}$ is a normal measure on $\check{\kappa}$.  I claim that for every $A \subset \kappa$ in $W$, either $q \forces \check{A}\in\dot{\mu}$ or $q \forces \check{A}\notin\dot{\mu}$.  If this is not the case, then there are $q_0,q_1 \leq q$ with $q_0 \Vdash \check{A}\in\dot{\mu}$ and $q_1 \Vdash \check{A}\notin\dot{\mu}$.  Let $g_0 \times g_1$ be $W$-generic for $Q\times Q$ such that $\<q_0,q_1>\in g_0 \times g_1$.  Then $W[g_0] \vDash A \in  (\dot{\mu})_{g_0}$ and $W[g_1] \vDash A \notin  (\dot{\mu})_{g_1}$.  Work in $W[g_0]$ and $W[g_1]$ to build $L[(\dot{\mu})_{g_0}]$ and $L[(\dot{\mu})_{g_1}]$ respectively.  The uniqueness of the model $L[\mu]$ implies that $L[(\dot{\mu})_{g_0}]=L[(\dot{\mu})_{g_1}]=L[\mu]$, and uniqueness of the normal measure in $L[\mu]$ implies that $(\dot{\mu})_{g_0}=(\dot{\mu})_{g_1}=\mu$. However, $A\in(\dot{\mu})_{g_0}$ but $A\notin (\dot{\mu})_{g_0}$, a contradiction. This shows that $W\cap \mu$ is definable in  $W$ as $\{A\subset \kappa \mid q \forces \check{A}\in\dot{\mu}\}$, and so $\kappa$ is measurable in $W$.  As $L[\mu]$ is the minimal model in which $\kappa$ is measurable, $L[\mu]\subset W$, and so $L[\mu] = W$.  Thus the forcing adding $h$ was trivial, and so $W \vDash \GA$.  The referee has pointed out that the proof can be generalized to show that whenever $W$ is an extender model satisfying ``I'm $K$'' and $K$ is absolute to forcing extensions, then $W$ satisfies the Ground Axiom.
\end{proof}

In addition, in many cases the core model $K$ satisfies the Ground Axiom.  However, certain other canonical models do not. For example, Ralf Schindler observed that methods of Woodin show that the least model $M_1$ of one Woodin cardinal is in fact a forcing extension of an inner model.   The similarities between the known models of the Ground Axiom suggest many questions about the consequences of the axiom.  As a starting point, we can consider the relationship of the Ground Axiom to the Generalized Continuum Hypothesis ($\GCH$) and to the assertion that every set is ordinal definable ($\VHOD$).

\begin{TestQuestion}\mbox{ }
\begin{enumerate}
\item \label{Ques:GA->GCH} Does the $\GA$ imply the $\GCH$?
\item \label{Ques:GA->HOD} Does the $\GA$ imply $\VHOD$?
\end{enumerate}
\end{TestQuestion}

Both questions have negative answers.  The former is answered below in Corollary \ref{C:GAand-GCH}, and the latter is answered in a subsequent paper \cite{HamkinsReitzWoodin:TheGroundAxiomandVnotHOD}, joint with Hamkins, Woodin, and myself.

Finally, the examples mentioned above leave open the question of compatibility of the Ground Axiom with various large cardinal hypotheses not covered by the Core Models, such as $\GA +$\emph{ there is a supercompact cardinal}.  In Section \ref{S:ForcingGA}, I will demonstrate  a method for building models of the Ground Axiom that accomodates large cardinals, giving relative consistency of the Ground Axiom with measurable and supercompact cardinals among others.  This method will also show the consistency of $\GA + \neg \GCH$.  In Section \ref{S:GAandGCH}, I consider an adaptation of the method to produce diverse models of $\GA + \GCH$.  In Section \ref{S:BedrockAxiom}, I will turn my attention to a related notion, the \emph{Bedrock Axiom} ($\BA$), which asserts that either the universe is a model of the Ground Axiom (a `bedrock model') or is a forcing extension of a such a model.  I will show the consistency of the negation of the Bedrock Axiom.

Throughout the following I will use blackboard bold $\P$ for proper class partial orders and standard text $P$ for set partial orders.  In function definitions, three dots $\from$ indicates that a function may be a partial function, e.g. $f\from \ORD \rightarrow \ORD$.  The partial order $Add(\gamma,\delta)$ consists of all functions $f\from \gamma\times\delta \rightarrow 2$ of size $<\gamma$, and $Coll(\gamma,\delta)$ consists of functions $f\from \gamma \rightarrow \delta$ of size $<\gamma$. In describing initial segments of models other than $V$, such as forcing extensions, I will use $V[G]_\alpha$ to denote $\(V_\alpha)^{V[G]}$.

\section{The Ground Axiom is first-order expressible}\label{S:GAExpressible}

I begin by showing that these notions are, indeed, first-order expressible.  

\begin{Thm}\label{T:GAIsFirstOrder}
There is a first-order formula which holds in a model of set theory exactly when that model is a forcing extension of an inner model by nontrivial set forcing.
\end{Thm}

Theorem \ref{T:GAIsFirstOrder} is a consequence of Theorem \ref{T:GADefinition}, which gives a more detailed result.  A similar result is implicit in independent work of Woodin \cite{Woodin:GenericMultiverse}.

The formula asserting the Ground Axiom will be given explicitly below, but some definitions are required.  The first are the \emph{$\delta$ cover and $\delta$ approximation properties}, formulated by Hamkins \cite{Hamkins2003:ExtensionsWithApproximationAndCoverProperties}, which provide a framework for analyzing  extensions and inner models.

\begin{Def}\rm\emph{(Hamkins)}.  Suppose that $W\subseteq V$ are transitive models of (some fragment of) $\ZFC$, and $\delta$ is a cardinal in $V$.
\begin{enumerate}
\item $\<W,V>$ has the \emph{$\delta$ cover property} if and only if for each $A\in V$ with $A\subset W$ and $|A|^{V} < \delta$ there is a covering set $B$ in $W$ with $A\subset B$ and $|B|^{W} < \delta$.
\item $\<W,V>$ has the \emph{$\delta$ approximation property} if and only if for each $A\in V$ with $A\subset W$, if $A\cap B\in W$ for every $B\in W$ with $|B|^{W} < \delta$, then $A\in W$.
\end{enumerate}
\end{Def}

As I will be working in initial segments of the universe, I will also need an appropriate variant of $\ZFC$.

\begin{Def}\rm Let $\ZFC_{\delta}$ be the theory consisting of Zermelo Set Theory, Choice,  and $\leq\!\delta$-Replacement (that is, Replacement holds for functions with domain $\delta$, a regular cardinal), together with the axiom $$(*) \quad \forall A \; \exists \alpha \in ORD \; \exists E\subseteq\alpha\times\alpha \thickspace \thickspace \<\alpha, E> \cong \<\tc(\{A\}), \in>$$ which asserts ``every set is coded by a set of ordinals.''
\end{Def}

Formally, $\ZFC_\delta$ is a theory in the language of set theory together with a symbol for $\delta$, and includes the assertion that $\delta$ is a regular cardinal.  In any model $V$ of $\ZFC$, if $\gamma$ is a $\beth$-fixed point of cofinality $>\!\delta$, then $V_\gamma \vDash \ZFC_\delta$.  The main result of this section is given in Theorem \ref{T:GADefinition}, providing an explicit first-order statement $\Phi(\delta,z,P,G)$ which holds if and only if the universe is a set forcing extension of an inner model (that is, if and only if $V\vDash \neg \GA$).

\begin{Thm}\label{T:GADefinition} The Ground Axiom is first-order expressible.\\
Specifically, the Ground Axiom fails if and only if there exist $\delta, z, P,$ and $G$ satisfying the following statement:\\
\begin{description}
	\item[$\Phi(\delta,z,P,G)$] \mbox{ }\\$\delta$ is a regular cardinal, $P\in z$ is a poset of size $<\delta$, $G$ is $z$-generic for $P$, and for every $\beth$-fixed point $\gamma > \delta$ of cofinality $>\delta$, there exists a transitive structure $\M$ of height $\gamma$ such that:
		\begin{enumerate}
		\parskip=0pt
		\item \label{T:GADefinition:ZFCdelta}$\M$ is a model of $\ZFC_{\delta}$,
		\item \label{T:GADefinition:z=Hdelta+}$z = (H_{\delta^+})^\M$,
		\item \label{T:GADefinition:forcingext}$\M[G] = V_{\gamma}$, and
		\item \label{T:GADefinition:CoverandApproximation}$\M \subset V_{\gamma}$ satisfies the $\delta$ cover and $\delta$ approximation properties.
	\end{enumerate}
\end{description}
\end{Thm}

This result is closely related to a result of Laver \cite{Laver:CertainVeryLargeCardinalsNotCreated}, appearing below as Theorem \ref{T:VDefinable}, which shows that every model of set theory is a definable class in all of its set forcing extensions (independently observed by Woodin \cite{Woodin:GenericMultiverse}).  The two directions will be stated and proven separately as Lemmas \ref{L:GADefinition1} and \ref{L:GADefinition2}.

\begin{Sublemma}\label{L:GADefinition1} If the Ground Axiom fails, then there exist $\delta, z, P,$ and $G$ satisfying $\Phi(\delta,z,P,G)$.
\end{Sublemma}

\begin{proof}
Suppose $V=W[G]$ is a forcing extension of $W$ by a poset $P\in W$. Let $\delta = (|P|^+)^V$ and let $z=(H_{\delta^+})^W$.  I will argue that $\Phi(\delta,z,P,G)$ holds.  Fix $\gamma > \delta$ a $\beth$-fixed point of cofinality $>\!\delta$.  I will show that $W_\gamma$ witnesses the properties of the structure $\M$ for $\gamma$.  For property \ref{T:GADefinition:ZFCdelta}, note that $W_\gamma \vDash \ZC$ (Zermelo set theory and the Axiom of Choice) for any limit ordinal $\gamma$.  As $\gamma$ is a $\beth$-fixed point, every set in $W_\gamma$ has transitive closure of size $<\gamma$ and so is coded as a set of ordinals in $W_\gamma$.  That $W_\gamma$ satisfies $\leq\!\delta$-Replacement follows from the cofinality of $\gamma$, for any function $f: a\to W_\gamma$ with $|a|=\delta$ has its range contained in some $W_\beta$ for $\beta < \gamma$, and so $\ran(f)\in W_{\beta+1}$.  Thus $W_\gamma \vDash \ZFC_\delta$.  That $z=(H_{\delta^+})^{W_\gamma}$ holds by definition, and so $W_\gamma$ witnesses property \ref{T:GADefinition:z=Hdelta+}.  To see that $W_{\gamma}[G] = V_\gamma$, note that any $x\in W[G]$ of rank $<\!\gamma$ will have a name $\dot{x}\in V$ of rank $<\!\gamma$, provided $\gamma$ is a limit ordinal greater than $\rank(P)$.

It remains to consider the $\delta$ cover and $\delta$ approximation properties.  I will first show that $W \subset V$ satisfies these properties, and will then argue that the properties are inherited by $W_\gamma \subset V_\gamma$.  The fact that any forcing extension by forcing of size less than $\delta$ will satisfy the $\delta$ cover and $\delta$ approximation properties is a special case of a more general result of Hamkins \cite{Hamkins2003:ExtensionsWithApproximationAndCoverProperties}, and the proof presented here follows his strategy. That $W\subset V$ satisfies the $\delta$ cover property is a well known result for any forcing with the $\delta$-c.c..  For the $\delta$ approximation property, it suffices to show that the property holds for sets of ordinals, since under $\ZFC$ every set has a bijection with an ordinal.  Suppose $A\in V$ is a subset of $\theta$ and for every $C\in W$ such that $|C|^W<\delta$, the intersection $C\cap A$ is in $W$.  If $\dot{A}\in W$ is a $P$-name for $A$, then there is $q\in G$ such that $q\Vdash $``$\dot{A}\subset\check{\theta} \wedge{}$ every $\delta$-approximation of $\dot{A}$ is in $W$.''  To see that the latter property is first-order expressible, note that it is equivalent to consider only those $C\in\mathcal{P}(\theta)^W$.  Suppose further that $A\notin W$, and so there is $p \in G$ extending $q$ such that $p\forces A\notin W$.  As $|P| < \delta$, there is an enumeration $\{p_\beta \mid \beta < \eta\}$ for some $\eta < \delta$ of all conditions below $p$.  Working in $W$, I can choose for each $p_\beta$ an $\alpha_\beta \in \theta$ and $p^0_\beta, p^1_\beta < p_\beta$ satisfying $p^0_\beta \Vdash \alpha_\beta \notin \dot{A}$ and $p^1_\beta \Vdash \alpha_\beta \in \dot{A}$.  Let $C = \{\alpha_\beta \mid \beta < \eta\}$.  Note $C\in W$ has size $<\delta$, so $C\cap A=B\in W$ and there is $q\in G$ forcing $\check{C}\cap\dot{A}=\check{B}$.  In particular, $q$ decides the statement ``$\alpha_\beta \in \dot{A}$'' for every $\beta < \eta$.  Without loss of generality $q \leq p$, for if this is not the case, then $G$ contains a common extension of $p$ and $q$.  Thus $q=p_\beta$ for some $\beta < \eta$.  By construction, there are $p^0_\beta, p^1_\beta < p_\beta$ satisfying $p^0_\beta \Vdash \alpha_\beta \notin \dot{A}$ and $p^1_\beta \Vdash \alpha_\beta \in \dot{A}$.  However $q$ decides the value of ``$\alpha_\beta \in \dot{A}$,'' contradiction.  Thus $A\in W$, and so $W\subset V$ satisfies the $\delta$-approximation property.  Finally, to see that these are inherited by $W_\gamma \subset V_\gamma$ is straightforward.  For example, if a set $A\in V_\beta$ for some $\beta < \gamma$ is covered by a set $B\in W$, then taking $B\cap W_\beta$ gives a cover of $A$ in $W_\gamma$.  A similar argument shows that for a set in $V_\gamma$, if every $W_\gamma$ approximation of that set over $V_\gamma$ lies in $W_\gamma$, then every $W$ approximation of that set over $V$ lies in $W$.  Thus $W_\gamma \subset V_\gamma$ satisfies the $\delta$ cover and $\delta$ approximation properties.  
\end{proof}

For the reverse direction, a key step in the argument is that for each $\gamma$ the witnessing structure $\M$ is unique.  This is an application of a result of Laver \cite{Laver:CertainVeryLargeCardinalsNotCreated}, who showed that the power set of $\delta$ together with the $\delta$ cover and $\delta$ approximation properties uniquely determine an inner model.

\begin{Sublemma}\emph{(Laver)}\label{L:M=M'}
Suppose $\M, \M^\prime$ and $V$ are transitive models of $\ZFC_{\delta}$ for $\delta$  a regular cardinal of $V$, both $\M$ and $\M^\prime$ are submodels of $V$ for which $\M\subseteq V$ and $\M^\prime \subseteq V$ satisfy the $\delta$ cover and $\delta$ approximation properties, $\M$ and $\M^\prime$ have the same power set of $\delta$,  and $(\delta^{+})^{\M} = (\delta^{+})^{\M^\prime} = (\delta^{+})^{V}$.  Then $\M=\M^\prime$.
\end{Sublemma}

\begin{proof} (Hamkins)
Since $\M$ and $\M^\prime$ satisfy axiom $(*)$, in order to conclude $\M=\M^\prime$ it suffices to show that $\M$ and $\M^\prime$ have the same sets of ordinals.  Observe that for any $A\subset ORD$, the statement ``$|A|<\delta$'' is unambiguous between the models $\M, \M^\prime$, and $V$, for by the $\delta$ cover property if it is true in one, then it will be true in all others containing $A$.  Also, since $(\delta^{+})^{\M} = (\delta^{+})^{\M^\prime} = (\delta^{+})^{V}$ the statement $|A|=\delta$ is unambiguous as well.

I will first establish a certain simultaneous cover property, namely, that for every set of ordinals $A\in V$ of size $<\delta$ there is a set $B\in \M\cap \M^\prime$ of size $\leq\delta$ with $A\subseteq B$.  Let $A\subset \alpha$ be such a set and fix in $V$ a well-ordering $\preccurlyeq$ of $\mathcal{P}(\alpha)^{V}$.  Construct in $V$ a sequence $\<B_\xi \mid \xi < \delta>$ of subsets of $\alpha$, each of size $<\delta$.  Let $B_0 = A$.  If $B_\xi \in \M$, then let $B_{\xi+1} \in \M^\prime$ be the $\preccurlyeq$-least subset of $\alpha$ such that $B_{\xi+1}\supset B_\xi$ and $|B_{\xi+1}|<\delta$.  The $\delta$ cover property for $\<\M^\prime,V>$ guarantees the existence of such a set.  If $B_\xi \notin \M$, then let $B_{\xi+1} \in \M$ to be the  $\preccurlyeq$-least subset of $\alpha$ with $B_{\xi+1}\supset B_\xi$, and $|B_{\xi+1}|<\delta$.  For limit $\xi$, let $B_\xi = \bigcup_{\beta<\xi} B_\beta$ and note that $B_\xi \in V$ by $\leq\delta$-Replacement.  As $\delta$ is regular and each $B_\beta$ has size $<\delta$, I have $|B_\xi|<\delta$ as well.  Thus the construction can continue through every stage $\xi < \delta$.  Note that $B_\xi\in \M$ for cofinally many $\xi$, and this is true for $\M^\prime$ as well.  In $V$, let $B = \bigcup_{\xi<\delta} B_\xi$.  To show $B\in \M$, I will show that every $\delta$ approximation of $B$ is in $\M$.  If $C\in \M$ has size $<\delta$, then by regularity of $\delta$ the intersection $C\cap B$ is equal to $C\cap B_\xi$ for some $\xi<\delta$.  Without loss of generality $B_\xi\in \M$ as this occurs cofinally in the sequence, and so $C\cap B_\xi \in \M$.  Thus $C\cap B \in \M$ for every set $C$ of size $<\delta$ in $\M$, and so by the $\delta$ approximation property $B\in \M$.  The same proof shows that $B\in \M^\prime$.  This establishes the desired simulataneous cover property.

I next claim that $\M$ and $\M^\prime$ have the same sets of ordinals of size $<\delta$.  Fix $A$ a set of ordinals in $\M$ of size $<\delta$.  The simultaneous cover property gives $B\in \M\cap \M^\prime$ of size $\leq\delta$ with $A\subseteq B$.  Since $B$ has order type smaller than $\delta^{+}$ there is a well-ordering $w\in \M$ of $\delta$, or possibly of a subset of $\delta$, of order type $ot(B)$.  Since  $w\subset \delta\times\delta$, the assumption that $\M$ and $\M^\prime$ have the same subsets of $\delta$ implies that $w$ must be in $\M^\prime$ as well.  The ordering $w$ induces an enumeration $\{b_\alpha \mid \alpha <\delta^\prime\}$ of $B$ for some $\delta^\prime\leq\delta$, and this enumeration is in $\M^\prime$ by $\leq\delta$-Replacement.  The set $\{\alpha <\delta^\prime \mid b_\alpha\in A\}$ is definable in $\M$, and since $\mathcal{P}(\delta)^\M=\mathcal{P}(\delta)^{\M^\prime}$ it exists in $\M^\prime$ as well.  As $A$ is definable from this set, $B$, and $w$, it must be the case that $A\in \M^\prime$.  The same argment shows that every set of ordinals in $\M^\prime$ of size $<\!\delta$  is also in $\M$.

Finally, it remains to show that $\M$ and $\M^\prime$ have the same sets of ordinals.  Fix $A$ a set of ordinals in $\M$.  I will show that every $\delta$ approximation of $A$ is in $\M^\prime$.  Fix $B \in \M^\prime$ of size $<\delta$.  The claim above shows that $B\in \M$ and so $A\cap B \in \M$, and applying the claim once more shows that $A\cap B \in \M^\prime$.   Thus every $\delta$ approximation of $A$ is in $\M^\prime$, and so $A\in \M^\prime$.  This shows $\M\subset \M^\prime$, and the reverse inclusion follows by the same argument, yielding $\M=\M^\prime$. 
\end{proof}

\begin{Sublemma}\label{L:GADefinition2} If there exist $\delta, z, P,$ and $G$ satisfying $\Phi(\delta,z,P,G)$, then
\begin{enumerate}
\parskip=0pt
\item \label{L:GADefinition2:uniqueness}For each $\gamma$ as in $\Phi$ there is a unique $\M$ witnessing properties (1) through (4) of $\Phi$, which I will denote $\M_{\gamma}$.
\item \label{L:GADefinition2:coherence}The $\M_\gamma$ form a coherent sequence: $\gamma < \gamma^\prime \to \left(\M_{\gamma^\prime}\right)_\gamma = \M_\gamma$. 
\item \label{L:GADefinition2:ZFC}If $\M = \bigcup_{\gamma} \M_{\gamma}$, then $\M  \vDash \ZFC$.
\item \label{L:GADefinition2:forcingext}$V = \M[G]$ is forcing extension of $\M$ by $P$.
\item \label{L:GADefinition2:GAfails}Consequently, the Ground Axiom fails in $V$.
\end{enumerate}
\end{Sublemma}

\begin{proof}
Fix $\gamma$ as in $\Phi$ and suppose $\M$, $\M^\prime$ are transitive structures of height $\gamma$ witnessing the properties stated in $\Phi$.  Note that $(\delta^{+})^{\M} = (\delta^{+})^{\M^\prime} = (\delta^{+})^{V}$, with the first equality holding by the equality of $H_{\delta^+}$ in $\M$ and $\M^\prime$, and the second equality holding by the fact that $V_\gamma=\M^\prime[G]$ is a forcing extension of $\M^\prime$ by forcing of size $<\!\delta$.  Thus Lemma \ref{L:M=M'} applies, and so $\M=\M^\prime$.  Thus the structure of height $\gamma$ witnessing the properties of $\Phi$ is unique, and will be denoted $\M_\gamma$.

To show the coherence of the $\M_\gamma$ sequence, suppose $\gamma < \gamma^\prime$ and $\M_{\gamma^\prime}$ is a witnessing structure for $\gamma^\prime$.  If $\mathcal{W}=(\M_{\gamma^\prime})_\gamma$, I must show that $\mathcal{W}$ is a witnessing structure for $\gamma$.  That $\mathcal{W}\vDash \ZC$ is straightforward, as $\mathcal{W}$ is a transitive initial segment of a model of $\ZC$ of height a limit ordinal.  The argument used in Lemma \ref{L:GADefinition1} shows that $\mathcal{W}\vDash \leq\!\delta$-Replacement, the pertinent property being $\leq\delta$-Replacement in the larger model $\M_{\gamma^\prime}$.  To see that every set in $\mathcal{W}$ is coded by a set of ordinals, it suffices to show that the transitive closure of every $x\in\mathcal{W}$ has size $<\gamma$ in $\M_{\gamma^\prime}$.  Fix such $x$, and suppose that there is an injection in $\M_{\gamma^\prime}$ from $\gamma$ to $\tc(x)$.  Then this injection is in $V$ as well, contradicting the assumption that $\gamma$ is a $\beth$-fixed point of $V$.  Thus $\mathcal{W}\vDash \ZFC_{\delta}$.  The $\delta$ cover and $\delta$ approximation properties are inherited by $\mathcal{W}\subset V_\gamma$ from $\M_{\gamma^\prime} \subset V_{\gamma^\prime}$ by the same argument given in the proof of Lemma \ref{L:GADefinition1}.  That $z = (H_{\delta^+})^{\M_{\gamma^\prime}} = (H_{\delta^+})^\mathcal{W}$ follows trivially from the definition of $\mathcal{W}$.  Finally, every $x\in V_\gamma$ has a $P$-name in $W_{\gamma^\prime}$ of rank less than $\gamma$, and so $\mathcal{W}[G] = V_\gamma$.  Thus $\mathcal{W}$ is a witnessing structure of height $\gamma$ for $\Phi$, and by uniqueness $\mathcal{W} = (\M_{\gamma^\prime})_\gamma = \M_\gamma$.  This establishes coherency of the $\M_\gamma$ sequence.

Note that $\M = \bigcup_{\gamma} \M_{\gamma}$ is a definable transitive class in $V$.  That $\M$ satisfies Extensionality, Pairing, Union, Power Set, Infinity, Regularity, and Choice can be checked in each instance by simply going to a large enough $\M_\gamma$.  The only difficulty lies in showing that $\M$ satisfies the  Replacement and Separation Schemes.   For Replacement it suffices to show Collection.  Fix $A\in\M$ and $\psi(x,y)$.  As $\M$ is definable in $V$ it follows that satisfaction of $\psi$ in $\M$ can be computed in $V$, that is, for $x,y$ in $\M$ I have $\M \vDash \psi(x,b)\leftrightarrow V \vDash \psi^\M(x,b)$. Since Replacement holds in $V$, there must be an $\alpha$ such that $\forall x\in A \thickspace \left(\exists y\in\M \thickspace \psi^\M(x,y) \to \exists y\in\M_\alpha \thickspace \psi^\M(x,y)\right)$.  Thus $\M \vDash \forall x\in A \thickspace  (\exists y \, \psi(x,y) \to \exists y\in\M_\alpha \, \psi(x,y))$, and so $\M$ satisfies Collection.  To see that $\M$ satisfies Separation, first note that a standard argument shows that $\M$ satisfies the reflection principle.  Fix $A\in\M$ and $\psi(x)$.  Choose $\alpha$ sufficiently large that $A \in \M_\alpha$ and $\M_\alpha$ reflects $\psi(x)$.  Choose $\gamma > \alpha$ a $\beth$-fixed point of cofinality $>\delta$.  Separation holds in $\M_\gamma$, so $B=\{x\in A \mid \M_\gamma \vDash \psi^{\M_\alpha}(x)\}$ is in $\M_\gamma$ and hence in $\M$.  But for $x\in A$, reflection gives $\M_\gamma \vDash \psi^{\M_\alpha}(x) \leftrightarrow \M_\alpha \vDash \psi(x) \leftrightarrow \M \vDash \psi(x)$.  Thus $B=\{x\in A \mid \M\vDash\psi(x)\}$, and so Separation holds in $\M$.  

It remains to show $V=\M[G]$.  Note that $G$ is $\M$-generic for $P$.  Furthermore, for any $A\in V$ there is find $\gamma$ above the rank of $A$ such that $A\in \M_\gamma[G]$, and so $A$ is the interpretation by $G$ of some $P$-name $\dot{A}\in\M_\gamma \subset \M$.  Similarly, for any $P$-name $\tau\in\M$ the valuation $\tau_G$ will be in $V$.  Thus $V=\M[G]$ is a forcing extension of $\M$, and so $V$ does not satisfy the Ground Axiom. 
\end{proof}

Note that if $V=\M[G]$ is a set forcing extension, Lemma \ref{L:GADefinition2} gives a first-order definition of $\M$ as a class of $V$ based on the parameters $\delta, (H_{\delta^+})^\M, P,$ and $G$.  By quantifying over $G$ we obtain Laver's result \cite{Laver:CertainVeryLargeCardinalsNotCreated}.

\begin{Thm}\label{T:VDefinable}\emph{(Laver)} 
Suppose $V=\M[G]$ is a forcing extension of $\M$ by set forcing $P \in \M$.  Then $\M$ is definable in $V$ from parameters in $\M$.  
\end{Thm}

\section{Forcing the Ground Axiom}\label{S:ForcingGA}

In order to explore the relative consistency of the Ground Axiom with other set theoretic assertions a general method for producing models of the Ground Axiom is required.  I will begin by identifying an axiom which implies the Ground Axiom.  Informally, this axiom asserts that every set in the universe is coded into the pattern of the $\GCH$ holding and failing at successor cardinals.  An exact statement of the assertion, which I will refer to as the Continuum Coding Axiom, is given below.  After demonstrating that the Ground Axiom is a consequence of the Continuum Coding Axiom, I will discuss a general method by which the property can be forced.   The apparent paradox -- using forcing to produce a model of the Ground Axiom -- is explained by noting that the Ground Axiom refers only to set forcing.  A proper class of forcing will be used to obtain the Continuum Coding Axiom.

\begin{Def}\rm\label{D:ContinuumCodingAxiom}The \emph{Continuum Coding Axiom} ($\CCA$) is the assertion that for every ordinal $\alpha$ and for every $a\subset \alpha$ there is an ordinal $\theta$ such that $\beta \in a \iff 2^{\aleph_{\theta+\beta+1}} = \(\aleph_{\theta+\beta+1})^+$ for every $\beta < \alpha$.
\end{Def}

Note that the $\CCA$ relies on coding only at the successor cardinals, over which we have greater control.  While the $\CCA$ refers only to sets of ordinals, this is equivalent in $\ZFC$ to every set being coded into the continuum function, since under $\AC$ every set is coded as a set of ordinals.  Furthermore, the $\CCA$ implies the apparently stronger assertion that every set of ordinals is coded \emph{cofinally often} into the power set function. Finally, note that the $\CCA$ is essentially a strong form of $\VHOD$.

\begin{Thm}\label{T:CCAimpliesGA} 
The $\CCA$ implies the $\GA$.
\end{Thm}

\begin{proof}
Suppose $V$ satisfies the Continuum Coding Axiom.  Suppose further that $V$ is a set forcing extension of an inner model $V=W[h]$, where $h$ is $W$-generic for some poset $Q\in W$. For  $\kappa > |Q|$, the models $W$ and $V$ will agree on the properties ``$\kappa$ is a cardinal'' and ``the $\GCH$ holds at $\kappa$.''  Every set of ordinals $a$ in $V$ is coded into the continuum function of $V$.  I claim that one such code for $a$ must appear above $|Q|$.  If $|Q|=\aleph_\delta$, consider the set of ordinals $a^\prime=\{\delta+\beta \mid \beta \in a\}$.  As $a^\prime$ is also coded into the continuum function, it is clear that the part of $a^\prime$ above $\delta$ must appear coded into the continuum function above $|Q|$.  Thus $a$ is coded into the continuum function of $V$ above $|Q|$, and so the code appears also in $W$.  Thus $a\in W$, and so every set of ordinals of $V$ is also in $W$.  This shows that $V=W$, and so the forcing $Q$ was trivial.  Thus $V \vDash \GA$. 
\end{proof}

To obtain a general method for forcing the Ground Axiom, I will present a proper class notion of forcing which gives $\CCA$ in the resulting extension.  This folkloric result is based on the work of Kenneth McAloon who presented the basic method in his paper on Ordinal Definability\cite{McAloon1971:ConsistencyResultsAboutOrdinalDefinability}.  See Jech \cite{Jech:SetTheory3rdEdition} and others for the basics of class forcing.

\begin{Thm}\label{T:CCA}
If $V$ satisfies $\ZFC$, then there is a forcing extension by class forcing which satisfies $\ZFC + \CCA$.
\end{Thm}

\begin{proof}The basic idea is to use forcing to code every set in the universe into the continuum function in the manner of the $\CCA$.  Coding a single set into the continuum function can be accomplished by applying Easton's celebrated result concerning powers of regular cardinals \cite{Easton70:PowersofRegularCardinals}.  Recall that a function $E\from \ORD \rightarrow \ORD$ is an \emph{Easton index function} if $E$ is nondecreasing, $\dom(E)$ is a set of regular cardinals,  and $E(\kappa)$ is a cardinal with $\cof\({E(\kappa)})>\kappa$ for each $\kappa \in \dom(E)$.  Associated with $E$ is a poset $Q(E)$, called Easton forcing, consisting of functions $q$ with domain $\subset \dom(E)$ and, for every regular $\lambda$, the domain of $q$ is bounded in $\lambda$, satisfying for each $\kappa$ in the domain, $q(\kappa)$ is a partial function from $E(\kappa)$ to $2$ of size $<\!\kappa$.  Conditions are ordered by extension on each coordinate and by enlarging the domain, with trivial extensions and enlargements allowed.  Easton's theorem tells us that, under mild $\GCH$ restrictions in the ground model, forcing with $Q(E)$ preserves cardinals and yields a model in which the power set of $\kappa$ has size $E(\kappa)$ for every $\kappa \in \dom(E)$.   What is more, for all $\lambda \notin \dom(E)$ the power set of $\lambda$ is `as small as possible,' i.e. the least cardinal $\alpha \geq \lambda^+$ such that $\cof(\alpha)>\lambda$ and $\alpha\geq E(\kappa)$ for all $\kappa \in \dom(E)\cap\lambda$.  The forcing relies only on a local $\GCH$ assumption, for if $\dom(E)\cup \ran(E)$ is contained in the closed interval $[\delta,\gamma]$ for regular cardinals $\delta$ and $\gamma$, then Easton's Theorem holds if $2^{<\!\delta}=\delta$ and  $2^\kappa = \kappa^+$ for every $\kappa \in [\delta,\gamma)$.  To code a set of ordinals $a\subset \alpha$ into the continuum function starting at  $\aleph_\theta$, define the Easton index function $E_a \from \aleph_{\theta+\alpha} \rightarrow \aleph_{\theta+\alpha}$ such that for all $\beta<\alpha$, the value $E_a \(\aleph_{\theta+\beta+1})$ is $\(\aleph_{\theta+\beta+1})^{+}$ for $\beta\in a$, and $\(\aleph_{\theta+\beta+1})^{++}$ otherwise.  Provided the $\GCH$ hypotheses described above hold on the interval $[\aleph_{\theta+1}, \aleph_{\theta+\alpha+1}]$, forcing with $Q(E_a)$ will code $a$ into the continuum function in exactly the manner of the $\CCA$.  Furthermore, the continuum function will be undisturbed outside this interval.

To prove the Theorem, I will encode every set into the continuum function. This will be accomplished by forcing with an $\ORD$-length iteration which, at stage $\xi$, uses Easton's result to encode the generics obtained at all previous stages.  Note that although the iteration only encodes generics, an easy density argument will show that in fact every set is encoded.  The coding at stage $\xi$ must take place on a ``clean'' interval, lying above the coding performed at all previous stages, to prevent later stages from destroying the coding at earlier stages.  This is accomplished by breaking the cardinals up into intervals, each of sufficient length to encode the generics added on all previous intervals.  The $\xi^{th}$ interval will be given by $[\aleph_{f(\xi)+1},\aleph_{f(\xi+1)})$, where $f(\xi)$ is defined recursively by $f(0)=0$, at successors $f(\xi+1)=\aleph_{f(\xi)}$, and for $\lambda$ a limit $f(\lambda)= \left(\sup_{\alpha<\lambda} f\(\alpha)\right) +1$.

The forcing is defined as follows.  Let $\P=\bigcup P_{\xi}$ be the $\ORD$-length iteration with Easton support, that is bounded support at inaccessible $\xi$ and full support otherwise, satisfying the following property.  For each $\xi$, the poset $P_{\xi+1}=P_\xi*\dot{Q}_\xi$, where the forcing $\dot{Q}_\xi$ is (a $P_\xi$-name for) the Easton forcing encoding the generic $G_\xi \subset P_\xi$ into the interval $[\aleph_{f(\xi)+1},\aleph_{f(\xi+1)})$, that is $P_{\xi} \forces \dot{Q}_\xi = Q(E_{\dot{G}_\xi})$.  Note that $\dot{G}_\xi$ is the canonical name for the generic subset of $P_\xi$.  Also observe that $G_\xi$ will not, strictly speaking, be a set of ordinals, but will rather be a sequence of sets $g_\eta$, each generic for the Easton forcing at stage $\eta$.  Each $g_\eta$ is itself a sequence  $g_\eta = \<g_\eta(\kappa)\mid \kappa \in \dom(g_\eta)>$, with each $g_\eta(\kappa)$ a binary sequence of length $\leq \kappa^{++}$.  By concatenating the $g_\eta(\kappa)$, I ``flatten'' $g_\eta$ into a single set $\corners{{g}_\eta} = \{\kappa^{++}+\beta \mid g_\eta(\kappa)(\beta)=1 \mbox{ for }\kappa \in \dom(g_\eta) \mbox{ and } \beta<\kappa\}$.  Note that $\corners{{g}_\eta} \subset \aleph_{f(\eta+1)}$.  The sequence of $\corners{{g}_\eta}$ are themselves concatenated to code $G_\xi$ as a single set of ordinals $\corners{{G}_\xi} = \{\aleph_{f(\eta+1)}+\beta \mid \beta \in \corners{{g}_\xi} \mbox{ for } \eta < \xi \mbox{ and } \beta < \aleph_{f(\eta+1)}\}$.  Here $\corners{{G}_\xi} \subset \aleph_{f(\xi)}$.

Now suppose $G$ is $V$-generic for $\P$.  I will assume that the $\GCH$ holds in $V$, as  if this is not the case, then it can be forced as an initial step of the proof.  That $V[G]$ satisfies $\ZFC$ follows from Easton support together with the increasing closure of the stage $\xi$ forcing as $\xi$ progresses through the ordinals.  For each $\xi$, let $G_\xi$ be the corresponding $V$-generic for $P_\xi$.  I claim that in $V[G]$, I can define $\corners{{G}_\xi}$ as $\{\beta<\aleph_{f(\xi)} \mid 2^{\aleph_{f(\xi)+\beta+1}}=(\aleph_{f(\xi)+\beta+1})^+\}$.  Factor $\P = P_\xi*\dot{Q}_\xi*\P_{tail}$. A standard argument establishes $|P_\xi|\leq \aleph_{f(\xi)+1}$, and it follows that the interval $[\aleph_{f(\xi)+1},\aleph_{f(\xi+1)+1}]$ satisfies the $\GCH$ requirements for Easton's Theorem in $V[G_\xi]$.  Therefore forcing with $(\dot{Q}_\xi)_{G_\xi}$ over $V[G_\xi]$ codes $\corners{{G}_\xi}$ into $[\aleph_{f(\xi)+1}, \aleph_{f(\xi+1)+1}]$ in the manner of the $\CCA$.  The tail forcing is $\leq \aleph_{f(\xi+1)}$-closed, so the coding is preserved in the full extension $V[G]$.  Thus in $V[G]$, every $\corners{{G}_\xi}$ is encoded in the manner of the $\CCA$.

It remains to show that every set of ordinals $a\in V[G]$ is encoded in the same fashion.  Fix $a\subset \alpha$ in $V[G]$ and $\xi$ such $\alpha<(\aleph_{f(\xi)})^+$.   Factor $\P = P_\xi*\dot{Q}_\xi*\P_{tail}$.  The closure of $\dot{Q}_\xi*\P_{tail}$ ensures that $a \in V[G_\xi]$.  Furthermore, the forcing $Q_\xi$ adds a subset to $(\aleph_{f(\xi)})^+$ which is Cohen generic over $V[G_\xi]$.  A standard density argument shows that every bounded subset of $(\aleph_{f(\xi)})^+$ appears as a block in the resulting Cohen generic.  That is, if $g\subset(\aleph_{f(\xi)})^+$ is $V[G_\xi]$-generic for $Add((\aleph_{f(\xi)})^+,1)$, then there is an ordinal $\theta$ such that $\beta \in a \iff \theta+\beta+1 \in g$ for every $\beta < \alpha$.  The method used to ``flatten'' $g_{\xi}$ to $\corners{{g}_{\xi}}$ and $G_{\xi+1}$ to $\corners{{G}_{\xi+1}}$ preserves each of the generics added on each coordinate of $Q_\xi$ as a contiguous block, and so $a$ appears as a block in $\corners{{G}_{\xi+1}}$.  The forcing at stage $\xi+1$ codes $\corners{{G}_{\xi+1}}\subset \aleph_{f(\xi+1)}$ into the continuum function, and so $a$ is coded into the continuum function as well.  Thus in $V[G_{\xi+1}]$ there is an ordinal $\theta$ such that $\beta \in a \iff 2^{\aleph_{\theta+\beta+1}} = (\aleph_{\theta+\beta+1})^+$ for every $\beta < \alpha$.  Thus $V[G] \vDash \CCA$.
\end{proof}

As the $\GCH$ fails cofinally in $V[G]$, an immediate consequence is

\begin{Cor}\label{C:GAand-GCH} If $V$ satisfies $\ZFC$, then there is a forcing extension by class forcing which satisfies $\ZFC + \GA + \neg \GCH$.  Thus, the answer to Test Question \ref{Ques:GA->GCH} is negative.
\end{Cor}

Furthermore, the method for constructing $V[G]$ can be adapted to provide additional consistency results.  For example, the encoding does not have to begin at $\aleph_0$. The iteration can be easily modified to begin at $\delta$ for any regular $\delta$, which means that the entire iteration will be $<\!\delta$-closed and so no sets will be added to $H_\delta$.  Thus an arbitrary initial segment of the universe can be preserved while forcing the Ground Axiom.

\begin{Thm}\label{T:GAandPreserveVAlpha}If $V \vDash \ZFC$ and $\alpha$ is an ordinal, then there is a forcing extension $V[G]$ of $V$ satisfying the Ground Axiom and having the same initial segment of height $\alpha$, i.e. $V[G]_\alpha = V_\alpha$.
\end{Thm}

\begin{proof}  Working in $V$, fix $\delta > |V_\alpha|$.  If the $\GCH$ does not hold above $\delta$, then as a first step it can be forced to hold there with $<\!\delta$-closed forcing using standard methods.  Then carry out the coding forcing described above, modifying the iteration to begin at $\delta$.  The resulting extension will satisfy the $\CCA$ and thus the Ground Axiom, and as both the $\GCH$ forcing and the coding forcing are $<\delta$-closed no sets will have been added to $V_\alpha$.
\end{proof}

This allows, for example, a measure on $\kappa$ to be preserved by the forcing, because a measure on $\kappa$ is verified in $V_{\kappa+2}$.

\begin{Cor}\label{C:GAand-GCHandMsbl} If $V$ satisfies $\ZFC + \kappa$ is measurable, then there is a forcing extension by class forcing which satisfies $\ZFC + \kappa$ is measurable $+ \GA + \neg \GCH$. 
\end{Cor}

This fact holds not just for `existence of a measurable cardinal' but for any $\Sigma_2$ property.

\begin{Cor}\label{C:GAisSigma2Independent1} If $\phi$ is any $\Sigma_2$ assertion true in $V$, then there is a forcing extension of $V$ by class forcing satisfying $\phi + \GA + \neg \GCH$.
\end{Cor}

\begin{proof}
This is a consequence of the following characterization of $\Sigma_2$ properties. Given a formula $\phi$, the following are equivalent.
\begin{enumerate}
\parskip=0pt
\item \label{C:Sigma2Equivalence1}$\ZFC \proves \phi \iff \sigma$, for some $\Sigma_2$ formula $\sigma$.
\item \label{C:Sigma2Equivalence2}$\ZFC \proves \phi \iff \exists \kappa \; H_\kappa \vDash \psi$, for some first-order formula $\psi$.
\end{enumerate}
The Corollary follows directly. Given a $\Sigma_2$ assertion $\exists x\forall y\: \phi_0(x,y)$ true in $V$, I choose $\kappa$ so that $H_\kappa \vDash \exists x\forall y\: \phi_0(x,y)$.  Then apply Theorem \ref{T:GAandPreserveVAlpha} to force the $\GA$ while preserving $H_\kappa$.  The resulting extension $V[G]$ will satisfy $\exists \kappa \, H_\kappa \vDash \exists x\forall y\: \phi_0(x,y)$, and so $V[G] \vDash  \exists x\forall y\: \phi_0(x,y)$.
\end{proof}

One suprising application of Corollary \ref{C:GAisSigma2Independent1} is the consistency of the Ground Axiom with Martin's Axiom.  Intuitively, Martin's Axiom seems to say ``a lot of forcing has been done,'' whereas the Ground Axiom seems to say the opposite.  In fact, an even stronger result along these lines holds, that the Ground Axiom is relatively consistent with the Proper Forcing Axiom.

\begin{Cor}\label{C:PFA+GA}If $V$ satisfies the Proper Forcing Axiom \emph{($\PFA$)}, then there is a class forcing extension of $V$ satisfying $\PFA + \GA$.
\end{Cor}

\begin{proof}This follows from the fact that the $\PFA$ is indestructible by $<\aleph_2$-directed closed forcing, (Larson \cite{Larson00:SeparatingStationaryReflectionPrinciples}).  By forcing the $\CCA$ as above but starting the iteration at $\aleph_2$, the iteration will be $<\aleph_2$-directed closed, preserving the $\PFA$ and forcing the $\GA$.
\end{proof}

Theorem \ref{T:GAandPreserveVAlpha} leaves open the question of consistency of the Ground Axiom with supercompact cardinals and other axioms not captured by $\Sigma_2$ formulas.  Just as in the previous corollary, however, the technology of indestructibility can be used to obtain further results.  Since we generally expect the Ground Axiom to hold in the canonical inner models of large cardinals, this theorem fits into the set theoretic program of obtaining large cardinal inner model properties by forcing.

\begin{Thm}\label{T:GAand-GCHandSupercompact} If $V$ satisfies $\ZFC + \kappa$ is supercompact, then there is a forcing extension by class forcing which satisfies $\ZFC + \kappa$ is supercompact $+ \GA + \neg \GCH$.
\end{Thm}

\begin{proof}
Laver's well-known result on indestructibility \cite{Laver78} shows that the supercompactness of $\kappa$ can be made indestructible by $<\!\kappa$-directed closed forcing.   For simplicity in the proof, I will assume that this condition holds in $V$.  Now force $\CCA$ as above with an iteration $\P$ encoding every set into the continuum function, but start the iteration at $\kappa$.  I will argue that $\kappa$ remains supercompact in the resulting extension $V[G]$.

Note that the entire iteration $\P$ is  $<\!\kappa$-directed closed.  However, Laver's result applies only to set forcing so an additional argument is needed.  Fix $\theta > \kappa$.  I will show that the $\theta$-supercompactness of $\kappa$ is preserved in $V[G]$.  Fix $\beta$ large enough so that $\P$ factors as $P_\beta * \P_{tail}$ where the second factor is $\leq \theta^{<\kappa}$-closed.  The first factor $P_\beta$ is $<\!\kappa$-directed closed set forcing, and so $\kappa$ remains $\theta$-supercompact in the partial extension $V[G_\beta]$.  The tail forcing adds no subsets to $\mathcal{P}_\kappa(\theta)$ and so $\kappa$ remains $\theta$-supercompact in the full extension $V[G]$.  Since this is true for every $\theta$, I have $\kappa$ supercompact in $V[G]$.  
\end{proof} 

This method generalizes to other large cardinals for which a comparable indestructibility theorem exists.

The methods for obtaining models of the Ground Axiom described above and in the following section work according to the same basic principle, forcing $\CCA$ or some similar coding axiom.  In particular, the models of the $\GA$ thus produced all satisfy strong versions of $\VHOD$.  It is natural to consider the relationship of $\VHOD$ to the $\GA$.  This relationship arose in Test Question \ref{Ques:GA->HOD}, which asked if the $\GA$ implies $\VHOD$.  The converse question, whether $\VHOD$ implies the $\GA$, is also natural.

In fact, neither implication holds.  Consistency of the $\GA$ with $\nVHOD$ is demonstrated in a subsequent paper \cite{HamkinsReitzWoodin:TheGroundAxiomandVnotHOD}, joint with Hamkins, Woodin, and myself, in which we show that every model of $\ZFC$ has forcing extension satisfying $\GA+\nVHOD$.  Consistency of $\neg\GA$ with $\VHOD$ was demonstrated by McAloon in 1970 \cite{McAloon1971:ConsistencyResultsAboutOrdinalDefinability}.  His result introduced the idea of coding information into the regular cardinals, and using this idea he produced a set forcing extension $L[G]$ of $L$ in which the generic $G$ is ordinal definable.  Since every element in $L[G]$ is definable from $G$ together with a name from $L$, every element of $L[G]$ is ordinal definable.  As $L[G]$ is an extension by set forcing, clearly $L[G] \vDash \VHOD+\neg\GA$.

Note that the above idea will work for any set forcing over $L$ for which the generic is definable in the extension.  For example, Fuchs and Hamkins in \cite{FuchsHamkins:DegreesOfRigidity} describe Suslin trees with the unique branch property.  Forcing with such a tree $T$ adds a single branch through the tree, and so the generic is definable as the unique branch of $T$ in the extension.  Such trees exist in $L$, and forcing over $L$ with the $L$-least such Suslin tree yields an extension in which the generic is definable without parameters.  The extension will satify $\neg\GA + \VHOD$.

In fact, it is not necessary to work in $L$.  These techniques can be combined with the method for forcing the Ground Axiom to work over any model $V$.  

\begin{Thm}\label{T:notGA+V=HOD2}
If $V\vDash \ZFC$, then there is a class forcing extension $V[G][H]$ satisfying  $\ZFC + \neg\GA + \VHOD$.
\end{Thm}

\begin{proof}
Using the methods of Theorem \ref{T:CCA}, go to a forcing extension $V[G]$ by class forcing in which every set is coded into the continuum function.  Working in $V[G]$, I follow McAloon's strategy described above, doing set forcing to obtain an extension $V[G][H]$ of $V[G]$ in which $H$ is coded into the continuum function.  Since $H$ was added by set forcing, $V[G]$ and $V[G][H]$ agree on the continuum function above the size of the forcing that added $H$.  Thus every set in $V[G]$ remains coded into the continuum function in $V[G][H]$, and so every set in $V[G]$ is ordinal definable in $V[G][H]$.  Since every set in $V[G][H]$ is definable from $H$ together with a name from $V[G]$, it follows that every set is ordinal definable in $V[G][H]$. This completes the proof.  
\end{proof}

\section{The Ground Axiom and the GCH}\label{S:GAandGCH}

The coding used above to force the Ground Axiom is quite flexible, but has the feature that in the resulting model the $\GCH$ fails quite strongly. In this section I will explore a different method of coding which preserves the $\GCH$.  This method relies on the same basic strategy of coding each set into the successor cardinals, differing only in the coding mechanism used.  Rather than coding according to the $\GCH$ holding or failing at the $\beta^{th}$ cardinal, the coding here will be accomplished by controlling whether the $\beta^{th}$ cardinal of the ground model is collapsed in the extension.  Unfortunately, this requires some method of computing in the extension the $\beta^{th}$ cardinal of the ground model.  This will limit the method to work only over models with a certain absoluteness property.

\begin{Def}\rm\label{D:AbsoluteInnerModel}
An inner model $U\subset V$ is \emph{ordinal definable} if it is defined by a first-order formula $\psi$ with ordinal parameters.  In addition, $U$ is \emph{forcing robust} if and only if for any model $N$ of set theory and any forcing extension $N[G]$ such that every set in $N[G]$ is set generic over N, $U^N = U^{N[G]}$.  
\end{Def}

Observe that forcing robustness applies to all set forcing extensions as well as a great many class forcing extensions, including iterations with progressively higher closure at each stage.  Note that the canonical models $L$, $L[0^\#]$, and $L[\mu]$ are all ordinal definable, forcing robust models.  

\begin{Thm}\label{T:GAandGCH1}
If $V=U[g]$ is a set forcing extension of $U$ an ordinal definable, forcing robust model, then there is a forcing extension of $V$ by nontrivial class forcing which satisfies $\GA + \GCH$.
\end{Thm}

\begin{proof}To demonstrate the basic method, I will assume first that the forcing adding $g$ was trivial, that is $V=U$ is itself an ordinal definable, forcing robust model, and second that the $\GCH$ holds in $V$.  I will then describe the modifications to the argument necessary to avoid these assumptions.  Suppose that $V$ is an ordinal definable, forcing robust model satisfying the $\GCH$.  I will use an iteration to collapse cardinals of $V$, coding every set into the pattern of ``$(\aleph_\beta)^V$ is a cardinal'' holding or failing in the extension.  This coding differs from that described in the previous section in that the forcing at stage $\beta$ does not code all generics added at previous stages, but rather codes a single `bit' of a single generic added at a previous stage.  Furthermore, to avoid complications arising from collapsing many cardinals in a row, coding will take place not at every successor cardinal but at every other successor cardinal.  Let $\kappa_\beta$ be the $\beta^{th}$ odd cardinal of $V$, that is $\kappa_\beta$ is the $\beta^{th}$ cardinal of the form $\aleph_{\lambda+n}$ where $\lambda$ is a limit ordinal and $n<\omega$ is odd.  Let $\P=\bigcup P_\beta$ be the $\ORD$-length iteration with Easton support, such that for each $\beta$, the poset $P_{\beta+1}=P_\beta*\dot{Q}_\beta$, where the forcing $\dot{Q}_\beta$ at stage $\beta$ is trivial if $\beta \in \corners{{G}_\beta}$ and equals $Coll(\kappa_\beta,\kappa_\beta^+)$ if $\beta \notin \corners{{G}_\beta}$. The set $\corners{{G}_\beta}$ is obtained by ``flattening'' the generic $G_\beta$ for $P_\beta$ into a single set of ordinals in the following way.  For each $\xi < \beta$, the generic $g_\xi$ added at stage $\xi$ can be encoded as a single subset $\corners{{g}_\xi}\subset \kappa_\beta^+$ in a canonical way using a pairing function.  For concreteness, I will use the absolute pairing function in which pairs are ordered first by maximum, then by first coordinate, and then by second coordinate.  The entire generic $G_\beta$ for $P_\beta$ is then flattened into a single set of ordinals $\corners{{G}_\beta} = \left\{ \aleph_{\kappa_\xi} + \alpha \mid \xi < \beta \mbox{ and } \alpha \in \corners{{g}_\xi}\right\}$.  Informally, a copy of $\corners{{g}_\xi}$ appears in $\corners{{G}_\beta}$ in the interval $[{ \aleph_{\kappa_\xi}, \aleph_{\kappa_\xi^+ }})$.  Note that for $\beta < \beta^\prime$, the larger generic $\corners{{G}_{\beta^\prime}}$ end-extends the smaller $\corners{{G}_\beta}$.   

Now suppose $G$ is $V$-generic for $\P$.  That $V[G] \vDash \ZFC$ follows from the fact that for any $\beta$ the partial order factors $\P=P_\beta*\P_{tail}$, where $P_\beta \forces \P_{tail}$ is $<\kappa_\beta$-closed.  The same fact implies that every set in $V[G]$ is set generic over $V$, so forcing robustness of $V$ applies to the extension $V\subset V[G]$.  Standard factoring arguments show that $\P$ preserves the $\GCH$ and collapses only those cardinals of the form $\delta=\kappa_\beta^+$ such that $\beta \notin \corners{{G}_\beta}$.
It follows that for any $\beta$, I can define $\corners{{G}_\beta}$ in $V[G]$ as $\{\xi < \kappa_\beta \mid (\aleph_{\xi}^+)^V \mbox{ is a cardinal} \}$.  The definition requires the extension to calculate $(\aleph_{\xi}^+)^V$, which relies on the forcing robustness of $V$.  Similarly, the generic $\corners{g_\xi}$ for an individual stage of forcing $\xi < \beta$ appears as a block in $\corners{{G}_\beta}$ and is therefore coded into an interval of the regular cardinals via the same coding.     I next show that every set of ordinals $a$ of $V[G]$ is definable from $g_\beta$ for $\beta$ sufficiently large.  Fix $a\subset \gamma$ in $V[G]$ and take $\beta > \gamma$ sufficiently large that $a \in V[G_\beta]$.  Without loss of generality I can assume that forcing at stage $\beta$ is nontrivial and so has the form $Coll(\delta, \delta^+)$ for some regular $\delta > \gamma$.  A density argument shows that if $g$ is generic for $Coll(\delta, \delta^+)$, then there is a $\theta < \delta$ such that $\alpha \in a \iff g(\theta+\alpha)=0$ for all $\alpha<\gamma$.  Thus for sufficiently large $\beta$ the set $a$ is definable from parameters $g_\beta$ and $\theta$.  That $V[G]$ satisfies the Ground Axiom follows directly.  Suppose to the contrary that $V[G] = W[h]$, a forcing extension of $W$ by a poset $Q$.  Forcing robustness of $V$ implies that $V\subset W$ and that $W$ computes $V$ correctly.  As $Q$ cannot collapse cardinals greater than $|Q|$, the models $V[G]$ and $W$ agree on the statement ``$(\kappa_\xi^+)^V$ is a cardinal'' for $\xi$ sufficiently large.  Thus for large enough $\beta$ the generic $g_\beta$ will be coded in $W$.  Since every set of ordinals in $V[G]$ is definable from $g_\beta$ for sufficiently large $\beta$, I have $V[G]\subset W$.  Thus $V[G]=W$ and the forcing $Q$ was trivial, so $V[G]$ satisfies the Ground Axiom.


This completes the proof under the additional simplifying assumptions on $V$.  Now suppose $V$ is an ordinal definable, forcing robust model but $V$ does not satisfy the $\GCH$.  It is possible that the coding as described above may fail in this case.  In particular, if $\gamma < \delta$ has power set $>\delta^{+}$, then $Coll(\delta,\delta^+)$ will collapse all cardinals in the interval $[\delta^+, 2^\gamma]$ to $\delta$.  To avoid this issue, I will begin by forcing the $\GCH$ over $V$ in the canonical way, forcing with a proper class iteration that adds a single Cohen subset to each regular cardinal.  The resulting model $\overline{V}$ satisfies the $\GCH$, but of course $\overline{V}$ may no longer be an ordinal definable, forcing robust model.  However, it is completely determined in $V$ whether a given cardinal is collapsed in $\overline{V}$.  That is, if $\Q$ is the canonical forcing of the $\GCH$, then either $1_\Q \forces \kappa$ is a cardinal or $1_\Q \forces \kappa$ is collapsed, for every cardinal $\kappa$.  Thus any model that can compute $V$ can compute ``the $\beta^{th}$ odd cardinal of $\overline{V}$.''  Thus the coding described above can be carried out over $\overline{V}$ to obtain $\overline{V}[G]$, and the proof that $\overline{V}[G]$ satisfies $\GA + \GCH$ goes through as before.

Finally, suppose $V=U[g]$ is a forcing extension of an ordinal definable, forcing robust model $U$ by set forcing $Q \in U$.  As set forcing cannot collapse cardinals above the size of the forcing, it follows that $U$ and $V$ have the same cardinals above $|Q|$.  Begin by forcing the $GCH$ to obtain $\overline{V}$.  Note that $Q$ has no effect on whether the canonical forcing of the $\GCH$ collapses cardinals, at least for cardinals larger than $|Q|$.  Thus, for $\kappa > |Q|$ a regular cardinal, $U$ and $U[g]$ agree on the statement ``the canonical forcing of the $\GCH$ collapses $\kappa$.'' Next, perform the collapsing coding as described above, but begin the iteration after $|Q|$.  The resulting model $\overline{V}[G]$ will be able to correctly calculate the $\beta^{th}$ odd cardinal of $\overline{V}$, at least for cardinals above $|Q|$.  The remainder of the proof follows as above.
\end{proof}

Note that this result applies to any set forcing extension of $L$, $L[0^\#]$, $L[\mu]$, and $K$ under many hypotheses.  The method can be adapted to yield slightly more general results.  Both the coding iteration and the canonical forcing of the $\GCH$ can begin at any regular $\delta$, allowing the preservation of an arbitrary initial segment of the universe.

\begin{Cor}\label{C:GAandGCHoverL} Suppose $U$ is an ordinal definable, forcing robust model.  If $\phi$ is any $\Sigma_2$ assertion forceable over $U$ by set forcing, then there is a forcing extension of $U$ satisfying $\phi + \GA + \GCH$ holds beyond some cardinal $\delta$.
\end{Cor}

\begin{proof}[Proof of Corollary \ref{C:GAandGCHoverL}]  Suppose $\phi$ is a $\Sigma_2$ assertion that holds in $U[g]$, a set forcing extension of $U$ by $Q$.  Fix $\delta$ so that $H_\delta \vDash \phi$ and $\delta > |Q|$, and force the Ground Axiom using the above coding but beginning both the canonical forcing of the $\GCH$ and the coding iteration at $\delta^+$.  The resulting model $V[G]$ will satisfy $\GA+\GCH$ holds above $\delta$, and the argument given in the proof of Corollary \ref{C:GAisSigma2Independent1} shows that preservation of $H_\delta$ implies that $V[G]$ will also satisfy $\phi$.
\end{proof}

\section{The Bedrock Axiom}\label{S:BedrockAxiom}

What are the models of $\ZFC$ familiar to the working set theorist?  There are canonical models, including $L$, $L[0^\#]$, $L[\mu]$, the model $K$ under various hypotheses, and many more.  These models are generally characterized by some notion of `minimality,' and in many cases satisfy the Ground Axiom.  On the other hand, many consistency results can be established by performing set forcing over the minimal models, such as the  consistency of $\neg${\sc CH}, of Martin's Axiom, and so on.  These models clearly do not satisfy the Ground Axiom, but they `sit above' a model of the Ground Axiom with only set forcing separating them.  Indeed, a common property of many models of set theory is that they are either models of the Ground Axiom or set forcing extensions of such models.  Does this hold in general? To investigate this question, I define the \emph{Bedrock Axiom}.  In the usual forcing paradigm one starts in a ground model $V$ and does forcing to obtain an extension $V[G]$.  I would like to shift perspective to the extension. If a model $V$ is a set forcing extension of some inner model $W$, then I will refer to $W$ as a ground model of $V$.  If a ground model of $V$ satisfies the Ground Axiom, then I will call it a \emph{bedrock model} of $V$, for if one descends through the ground models of $V$ one reaches bedrock when one can descend no further.  The Bedrock Axiom asserts the existence of a bedrock model.

\begin{Def}\rm\label{D:BA}
The \emph{Bedrock Axiom} ($\BA$) asserts there is an inner model $W$ such that $V$ is a set forcing extension of $W$ and $W\vDash \GA$.
\end{Def}

In particular, any model of the Ground Axiom is a model of the Bedrock Axiom, since every model is trivially a forcing extension of itself.  Once again this apparently second-order statement has a first-order equivalent.

\begin{Thm}\label{T:BAExpressible}
The Bedrock Axiom is first-order expressible.
\end{Thm}

\begin{proof}
The Bedrock Axiom is expressed by the statement ``either $V$ satisfies the Ground Axiom or there exist $\delta, z, P$, and $G$ that satisfy the formula $\Phi$ of Theorem \ref{T:GADefinition} and, in the resulting inner model $W$, there are no $\delta^\prime, z^\prime, P^\prime$, and $G^\prime$ that satisfy $\Phi$ relativized to $W$.''
\end{proof}

Despite the many natural examples of models satisfying the Bedrock Axiom, it is consistent that the axiom fails.

\begin{Thm}\label{T:Con-BA}
There is a forcing extension of $L$ by class forcing which satisfies $\ZFC + \neg \BA$.
\end{Thm}

\begin{proof}
I will start in $L$ and build a class forcing extension, this time using a product rather than an iteration.  It is the commutative property of products that will be key in showing that the resulting model $L[G]$ satisfies $\neg \BA$.

For each regular cardinal $\lambda$ of $L$, let $P_\lambda = Add(\lambda,1)$.  Let $\P=\prod_{\lambda}P_\lambda$ be the canonical Easton product adding a single subset to each regular $\lambda$.  Easton's forcing is explored in detail in the standard texts (such as Jech \cite{Jech:SetTheory3rdEdition} and Kunen \cite{Kunen:Independence}), and I will use the basic results about the forcing without proof.  If $G$ is $V$-generic for $\P$, then I claim that $L[G]\vDash \neg \BA$.  I must show that if $L[G]$ is a forcing extension of some model $W$, then $W$ does not satisfy the Ground Axiom.  As a warm-up, I observe that $L[G]$ itself does not satisfy the Ground Axiom.  Note that for any regular $\lambda$, I can factor $\P$ as $\P\cong \P^{>\!\lambda}\times P^{\leq\!\lambda}$ where $P^{\leq\lambda}$ has the $\lambda^+$-c.c. and $\P^{>\!\lambda}$ is $\leq\!\lambda$-closed.  Thus $L[G]=L[G^{>\lambda}][G^{\leq\lambda}]$, a nontrivial set forcing extension of $L[G^{>\lambda}]$, and so $L[G]$ does not satisfy the Ground Axiom. 

Now suppose $L[G] = W[h]$ where $h$ is $W$-generic for some poset $Q\in W$.  Let $\delta = \left(|Q|^{+}\right)^W$.  Then $\delta$ is a regular cardinal of $L$, and $\P\cong \P^{>\delta}\times P^{\leq\delta}$.  I claim that $L[G^{>\delta}]\subset W$.  Clearly $L\subset W$, so it suffices to show that $g_\lambda\in W$ for regular $\lambda>\delta$.  Fix such $\lambda$ and consider $g_\lambda\subset \lambda$.  As shown in the proof of Lemma \ref{L:GADefinition1}, if $\lambda$ is a regular cardinal larger than $|Q|$, then $W \subset W[h]$ satisfies the $\lambda$ cover and $\lambda$ approximation properties.  Since $g_\lambda$ is $L$-generic for $Add(\lambda, 1)$, every initial segment of $g_\lambda$ is in $L$ and thus in $W$.  It follows that $g_\lambda\in W$ by the $\lambda$ approximation property, so $L[G^{>\delta}]\subset W$.

Thus $L[G^{>\delta}]\subset W\subset L[G]$, so $W$ is an intermediate model between  $L[G^{>\delta}]$ and a forcing extension $L[G^{>\delta}][G^{\leq\delta}]=L[G]$ by set forcing $P^{\leq\delta}$.  It follows that $W$ is a possibly trivial forcing extension of $L[G^{>\delta}]$, a result which appears in many places in the literature. For example in Jech \cite{Jech:SetTheory3rdEdition} p.265 it is stated as follows:

\begin{Sublemma}\label{L:FactoringBAs}
Let $G$ be $V$-generic on a complete Boolean algebra $B$.  If $N$ is a model of $\ZFC$ such that $V\subseteq N\subseteq V[G]$, then there exists a complete subalgebra $D\subseteq B$ such that $N=V[D\cap G]$. 
\end{Sublemma}

Continuing with the proof of Theorem \ref{T:Con-BA}, if $W\ne L[G^{>\delta}]$, then $W$ is a set forcing extension of $L[G^{>\delta}]$, and so does not satisfy the Ground Axiom.  If $W = L[G^{>\delta}]$, then observe that for any regular $\lambda>\delta$, I can factor $\P^{>\delta}\cong\P^{>\gamma}\times P^{(\delta,\gamma]}$, and $L[G^{>\delta}]=L[G^{>\gamma}][G^{(\delta,\gamma]}]$.  In this case, $W$ is a set forcing extension of $L[G^{>\gamma}]$ by set forcing $P^{(\delta,\gamma]}$ and so does not satisfy the Ground Axiom.  This shows no ground model $W$ of $L[G]$ satisfies the Ground Axiom, and so $L[G] \vDash \neg \BA$. 
\end{proof}

Note that in the above construction forcing need not occur at every regular $\lambda$ but may be restricted to a definable proper subclass.  The same result may be obtained by starting the forcing above any fixed $\kappa$, or by forcing at every other regular cardinal, etc.  This allows the technique to be combined with that used in Section \ref{S:ForcingGA} to provide results of much greater generality.

\begin{Thm}\label{T:-BA1}
If $V$ satisfies $\ZFC$, then there is a forcing extension by class forcing which satisfies $\ZFC + \neg \BA$.
\end{Thm}

\begin{proof}
Begin by following the strategy described in Section \ref{S:ForcingGA} to obtain a model $V[G]\vDash \CCA$ in which every set is definable from the continuum function.  Now proceed as in the proof above, forcing over $V[G]$ with a class product to add a single subset to regular cardinals $\lambda$.  However, in order to insure preservation of cardinals and of the continuum function we must force only at those regular $\lambda$ for which $2^{<\lambda}=\lambda$.  If $H$ is $V[G]$-generic for this product, then the resulting model $V[G][H]$ has the same cardinals and the same continuum function as $V[G]$. Thus in $V[G][H]$ every set in $V[G]$ is definable from the continuum function.  Now suppose $V[G][H]=W[h]$ where $h$ is $W$-generic for some poset $Q\in W$.  Once again, $W$ and $V[G][H]$ must agree on the value of $2^{\alpha}$ for $\alpha$ sufficiently large, as well as on the statement ``$\alpha$ is the $\beta^{th}$ cardinal.''  Thus every set in $V[G]$ is definable from the continuum function in $W$, and so $V[G] \subset W$.  Now set $\delta = (|Q|^{+})^W$.  For $\lambda > \delta$, I will show $h_\lambda \in W$.  As $V[G] \subset W$, every initial segment and hence every $\lambda$ approximation of $h_\lambda$ is in $W$.  Since $\<W,W[h]>$ satisfies the $\lambda$ approximation property, $h_\lambda \in W$.  Thus $V[G][H^{>\delta}]\subset W$.  I now have $V[G][H^{>\delta}] \subseteq W \subset V[G][H]$ and so, by factoring $H^{>\delta}$ if necessary, $W$ is a forcing extension of an inner model.  Thus $V[G][H] \vDash \neg \BA$.   
\end{proof}

Various modifications provide further results.  An arbitrary initial segment of $V$ can be preserved when constructing $V[G][H]$ by restricting both the class iteration and the class product to stages above some suitably chosen $\kappa$.  This allows the preservation of a measurable cardinal, as in Corollary \ref{C:GAand-GCHandMsbl}, and in combination with indestructibility allows preservation of a supercompact cardinal, as in Theorem \ref{T:GAand-GCHandSupercompact}. In addition, experts may find it natural to consider whether $M1$, the least iterable extender model with one Woodin cardinal, is a model of $\BA$.

\begin{Cor}\label{C:-BAandMsbl} Suppose $V$ satisfies $\ZFC$.  Then
\begin{enumerate}
\parskip=0pt
\item If $V$ satisfies $\kappa$ is measurable, then there is a forcing extension by class forcing which satisfies $\ZFC + \kappa$ is measurable  $+ \, \neg\BA$.
\item If $V$ satisfies $\kappa$ is supercompact, then there is a forcing extension by class forcing which satisfies $\ZFC + \kappa$ is supercompact $+ \, \neg\BA$. 
\end{enumerate}
\end{Cor}

In contrast to the models of the Ground Axiom, none of the models of $\neg\BA$ produced above satisfy $\VHOD$.   The forcing to add a subset to a cardinal is almost homogeneous (a poset is \emph{almost homogeneous} if for any $p,q$ there is a automorphism of the poset sending $p$ to condition compatible with $q$), and it follows that the full Easton product is as well.  It is known that extensions by almost homogeneous forcing always satisfy $V \neq \HOD$ (Kunen \cite{Kunen:Independence} pp.244-245).  However, by combining the product forcing above with a set version of the coding used in Theorem \ref{T:CCA} to force $\CCA$,  a model of $\VHOD + \neg\BA$ can be obtained.

\begin{Thm}\label{T:-BA+VHOD}There is a forcing extension of $L$ by class forcing satisfying $\VHOD + \neg\BA$.
\end{Thm}

\begin{proof}I will force with a class product $\P=\prod_{\lambda} P_\lambda$ with Easton support, that is bounded support at inaccessibles and full support otherwise.  Rather than each factor simply adding a set to a cardinal, the factor $P_\lambda$ will be a short iteration coding its own generic into the continuum function on a particular interval of cardinals.  To keep track of the intervals, I will once again use the function $f$ defined in Theorem \ref{T:CCA}, with $f(0)=0$, at successors $f(\xi+1)=\aleph_{f(\xi)}$, and for $\lambda$ a limit $f(\lambda)= \left(\sup_{\alpha<\lambda} f\(\alpha)\right) +1$.  Forcing will occur only at limit ordinals $\lambda$, and for such a $\lambda$ the forcing $P_\lambda$ is defined as $P_\lambda = \bigcup_{n\in\omega} P_{\lambda,n}$ a forcing iteration of length $\omega$ with full support.  For each $n<\omega$, the forcing $P_{\lambda,n+1}=P_{\lambda,n}*\dot{Q}_{\lambda,n}$ such that $\dot{Q}_{\lambda,n}$ is (a $P_{\lambda,n}$-name for) the Easton forcing coding the generic $G_{\lambda,n} \subset P_{\lambda,n}$ into the interval $[\aleph_{f(\lambda+n)+1},\aleph_{f(\lambda+n+1)})$.  Note that each generic $G_{\lambda,n}$ is ``flattened'' into a single set of ordinals $\corners{{G}_{\lambda,n}}$ before being encoded, just as in Theorem \ref{T:CCA}.  The generic $G_\lambda$ for $P_\lambda$ is completely determined by the collection of $\{G_{\lambda,n} \mid n\in \omega\}$, and forcing with $P_\lambda$ codes all $G_{\lambda,n}$ into the continuum function in the interval $[\aleph_{f(\lambda)+1},\aleph_{f(\lambda+\omega)})$ in the manner of the $\CCA$.  Furthermore, a standard analysis shows that $P_\lambda$ preserves cardinals and leaves the continuum function undisturbed outside of this interval.  The size and chain condition of $P_\lambda$ are both $\aleph_{f(\lambda+\omega)}$. For any limit $\lambda$, the product $\P$ factors as $\P \cong P^{\leq\lambda}\times \P^{>\lambda}$, where $P^{\leq\lambda}$ consists of all conditions $p\in\P$ such that $\dom(p)\subset\lambda+1$, and $\P^{>\lambda}$ consists of those conditions for which $\dom(p)\cap\lambda+1 = \emptyset$.  A straightforward calculation shows $P^{\leq\lambda}$ has size $\aleph_{f(\lambda+\omega)}$ and the $\aleph_{f(\lambda+\omega)}$-c.c., and $P^{>\lambda}$ is $\leq \aleph_{f(\lambda+\omega)}$-closed.  It is a standard fact from the theory of products that for $\kappa$ a regular cardinal, if $\P$ factors at $\kappa$, that is $\P\cong P^1\times \P^2$ where $P^1$ has the $\kappa^+$-c.c. and $\P^2$ is $\leq\kappa$-closed, then every function $f\from \kappa \rightarrow L$ in the extension by $\P$ is already in the extension by $P^1$ (Jech \cite{Jech:SetTheory3rdEdition} p.234).  Recalling that $\aleph_{f(\lambda+\omega)}$ is a successor cardinal and therefore regular, this gives arbitrarily large regular $\kappa$ such that $\P$ factors at $\kappa$.  This is sufficient to ensure that forcing with $\P$ preserves $\ZFC$.

Suppose $G$ is $L$-generic for $\P$.  I will first show that $\P$ preserves cardinals.  Fix $\kappa$ a regular cardinal of $L$.  If there exists a limit ordinal $\lambda$ such that $\aleph_{f(\lambda)} \leq \kappa < \aleph_{f(\lambda+\omega)}$, then factor $\P \cong P^{\leq\lambda}\times \P^{>\lambda}$.  If $\kappa$ is collapsed in $L[G]$, then it must already be collapsed in $L[G^{\leq\lambda}]$.  I now factor $P^{\leq\lambda} \cong P_\lambda \times P^{<\lambda}$, and so $L[G^{\leq\lambda}] = L[G_\lambda][G^{<\lambda}]$.  The forcing $P^{<\lambda}$ has size $\aleph_{f(\lambda)}$ and the $\aleph_{f(\lambda)}$-c.c. in $L$, and the closure of $P_\lambda$ means that this remains true in $L[G_\lambda]$.  Thus $P^{<\lambda}$ will not collapse $\kappa$ over $L[G_\lambda]$, so if $\kappa$ is collapsed it must already be collapsed in $L[G_\lambda]$.  This is impossible as $P_\lambda$ preserves cardinals, so $\kappa$ must remain a cardinal in $L[G]$.  If there is no limit $\lambda$ as described above, then $\kappa$ must be a limit of $\aleph_{f(\beta)}$ for $\beta<\lambda$.  If $\kappa$ is collapsed, then some successor $\delta<\kappa$ is also collapsed, which is impossible.  Thus $\P$ preserves all cardinals.

I claim that for every limit ordinal $\lambda$ and $n<\omega$, the code for the generic $\corners{{G}_{\lambda,n}}$ is definable in $L[G]$ as $\corners{{G}_{\lambda,n}} = \{ \beta < \aleph_{f(\lambda+n)} \mid  2^{\aleph_{f(\lambda+n)+\beta+1}} = ( \aleph_{f(\lambda+n)+\beta+1} )^+ \}$.  Fix $\lambda$ a limit and factor $\P \cong P^{\leq\lambda}\times \P^{>\lambda}$.  $\P^{>\lambda}$ is $<\aleph_{f(\lambda+\omega)}$-closed and so it suffices to show that the claim holds in $L[G^{\leq\lambda}]$.  I once again factor $P^{\leq\lambda} \cong P_\lambda \times P^{<\lambda}$, and as $P^{<\lambda}$ is too small to affect the $\GCH$ above $\aleph_{f(\lambda)+1}$ it suffices to show the claim holds in $L[G_\lambda]$.  However, the forcing $P_\lambda$ codes $\corners{{G}_{\lambda,n}}$ into the continuum function in exactly the way described in the claim.  Thus the claim holds in $L[G]$.  This shows that in $L[G]$ every $\corners{{G}_{\lambda,n}}$, and therefore every $G_{\lambda,n}$, is definable from parameters $\lambda$ and $n$.  As the definition is uniform, it follows that $G^{\leq\lambda}$ is ordinal definable in $L[G]$ for every $\lambda$.  Any set $a\in L[G]$ is in $L[G^{\leq\lambda}]$ for some $\lambda$, and so $a$ is definable from $G^{\leq\lambda}$ together with a name $\dot{a}\in L$ for $a$.  As every member of $L$ is ordinal definable in $L[G]$, $a$ is ordinal definable there as well.  Thus $L[G] \vDash \VHOD$.

It remains to show that $L[G]$ satisfies $\neg\BA$.  That $L[G]$ itself does not satisfy the Ground Axiom follows from the observation that for any $\lambda$ the model $L[G]$ can be written $L[G]=L[G^{>\lambda}][G^{\leq\lambda}]$, where the latter extension is by the set forcing $P^{\leq\lambda}$.  Now suppose $L[G]=W[h]$, a forcing extension of $W$ by a poset $Q$.  As $Q$ cannot affect the continuum function above $|Q|$, the models $W$ and $L[G]$ agree on the assertion $2^\delta = \delta^+$ for $\delta>|Q|$.  Thus for $\lambda$ such that $\aleph_{f(\lambda)} > |Q|$, the generic $G_\lambda$ is definable in $W$.  Fixing such a $\lambda$, I have $L[G^{>\lambda}]\subset W$.  The argument given in the proof of Theorem \ref{T:CCA} shows that $W$ is therefore a set forcing extension of $L[G^{>\lambda^\prime}]$ for $\lambda^\prime>\lambda$.  Thus $W$ does not satisfy the Ground Axiom, and so $L[G]\vDash \neg\BA$.
\end{proof}

\section{Open Questions}

A number of questions regarding these axioms remain.  Of particular interest is the question of uniqueness of bedrock models.

\begin{Question}If the Bedrock Axiom holds, is the bedrock model unique?
\end{Question}

A negative answer would settle another natural question:  Given two ground models of $V$ can we always find a third ground model contained in their intersection?  If the answer is yes, it would indicate that forcing can be used to amalgamate only models that are in some sense ``close together.''  Another natural structure to consider is the class obtained by intersecting all ground models.  Is it a model of $\ZFC$?  In many examples the answer is yes, even in cases as in the previous section where the Bedrock Axiom fails.  The models of $\neg\BA$ described above also have the property that there are a large number of ground models, proper class many.  Is it possible to have $\neg\BA$ with only set many ground models?  If so, can we reduce the number of ground models to countable?

Another area of interest is analyzing the restriction of these axioms to various classes of forcing such as {\sc ccc} or proper forcing.  For example, $\GA_{ccc}$ is the assertion that the universe is not a set forcing extension of an inner model by {\sc ccc} forcing.  Can we separate notions of the Ground Axiom by obtaining models of, for example, $\GA_{ccc} + \neg\GA_{proper}$?  How about consistency of restricted notions of $\GA$ and $\neg\BA$, such as $\GA_{ccc} + \neg\BA$?

Also of interest is the extension of these axioms to include class forcing.  While a first-order expression of such an axiom may be impossible in the general case, there is some hope that by restricting attention to a particular class of class forcing, for example, forcing with a closure point at $\delta$, a first-order expression could be achieved.

\bibliographystyle{alpha}
\bibliography{MathBiblio,HamkinsBiblio}

\end{document}